\documentclass[11pt]{article}
\usepackage[T1]{fontenc}
\usepackage{amsmath,amsthm,amssymb,dsfont,graphicx,xspace,epsfig,xcolor}
\usepackage{fullpage}
\usepackage{thm-restate}
\usepackage[hidelinks]{hyperref}
\usepackage{bm}
\usepackage{color}
\usepackage{comment}
\usepackage{mathrsfs}
\usepackage[shortlabels,inline]{enumitem}
\usepackage{float}
\usepackage{dutchcal}
\usepackage{pict2e}
\usepackage{mathtools}
\usepackage{tikz}
\usepackage{subfigure}
\usepackage{authblk}

\usetikzlibrary{fit, backgrounds, matrix, arrows.meta}
\usetikzlibrary{calc,arrows,positioning}
\usetikzlibrary{decorations}
\usetikzlibrary{decorations.pathmorphing}
\usetikzlibrary{decorations.pathreplacing}
\usetikzlibrary{decorations.shapes}
\usetikzlibrary{decorations.text}
\usetikzlibrary{decorations.markings}
\usetikzlibrary{decorations.fractals}
\usetikzlibrary{decorations.footprints}

\tikzset{vertex/.style = {circle,fill=black,minimum size=4pt, inner sep=0pt}}
\tikzset{smallvertex/.style = {circle,white,minimum size=1pt, inner sep=0pt}}
\tikzset{arc/.style = {thin,->,> = latex}}

\newtheorem{theorem}{Theorem}[section]
\newtheorem{corollary}[theorem]{Corollary}
\newtheorem{proposition}[theorem]{Proposition}
\newtheorem{lemma}[theorem]{Lemma}
\newtheorem{conjecture}[theorem]{Conjecture}
\theoremstyle{definition}
\newtheorem{problem}[theorem]{Problem}

\renewcommand{\epsilon}{\varepsilon}

\renewcommand{\phi}{\varphi}

\let\leq\leqslant
\let\geq\geqslant

\newcommand{\EE}{\mathbb{E}}
\newcommand{\PP}{\mathbb{P}}

\newcommand{\Geh}[1]{\mathcal{EH}^{(#1)}}
\newcommand{\directedK}{\vec{K}}
\newcommand{\directedKss}{\directedK_{s,s}}
\newcommand{\AP}{\mathcal{AP}}
\newcommand{\AC}{\mathcal{AC}}

\DeclareMathOperator{\Ad}{Ad}
\newcommand{\Ram}{R}

\newcommand{\fmain}{f_{\ref{thm:main}}}
\newcommand{\fmainoriented}{f_{\ref{thm:mainoriented}}}
\newcommand{\fKO}{f_{\ref{thm:KO}}}

\newcommand{\fKLST}{f_{\ref{thm:KLST}}}

\newcommand{\fcorsubstructureoriented}{f_{\ref{cor:substructure:oriented}}}

\hypersetup{
    colorlinks,
    linkcolor={red!50!black},
    citecolor={blue!50!black},
    urlcolor={blue!80!black}
}

\title{Edge-colouring and orientations: applications to degree- and $\chi$-boundedness\footnote{
		Research supported by JSPS KAKENHI JP20A402 and 22H05001, and by JST ASPIRE JPMJAP2302.}}

\author{Arnab Char}
\author{Ken-ichi Kawarabayashi}
\author{Lucas Picasarri-Arrieta\footnote{Corresponding author. Email: {\tt lpicasarr@nii.ac.jp}.}}

\affil{
National Institute of Informatics, The University of Tokyo, Tokyo, Japan. 
}

\date{}

\begin{document}

\maketitle

\begin{abstract}
We prove a new generalisation of Ramsey's theorem by showing that every $2$-edge-coloured graph with sufficiently large minimum degree contains a monochromatic induced subgraph whose minimum degree remains large. From this, we also derive that every orientation of a graph with large minimum degree contains either a large transitive tournament or an induced antidirected digraph whose minimum degree is still large.

As a consequence, we obtain two general tools showing that certain extensions of degree-bounded graph classes preserve degree-boundedness.
A hereditary class $\mathcal{G}$ is {\it degree-bounded} if, for every integer $s$, there exists $d=d(s)$ such that every graph $G\in \mathcal{G}$ either contains $K_{s,s}$ or has minimum degree at most $d$.
With these tools, we obtain for instance that odd-signable graphs and Burling graphs are degree-bounded. We also characterise exactly the oriented graphs $F$ such that the graphs admitting an orientation without any induced copy of $F$ are degree-bounded. 
\end{abstract}

\section{Introduction}

In this paper, we provide a new generalisation of Ramsey's theorem. Informally, we show that, in any $2$-edge-coloured graph with sufficiently large minimum degree, one can extract a monochromatic induced subgraph whose minimum degree remains large. 
Formally, a {\it $k$-edge-coloured} graph $G$ is a graph given with a $k$-colouring $\phi$ of its edges. A subgraph of $G$ is {\it monochromatic} if all its edges receive the same colour under $\phi$.

\begin{restatable}{theorem}{mainthm}
    \label{thm:main}
    There exists $\fmain\colon\mathbb{N}^2\to \mathbb{N}$ such that, for all $k,d\in \mathbb{N}$, every $k$-edge-coloured graph $G$ with $\delta(G) \geq \fmain(k,d)$ contains a monochromatic induced subgraph $H$ with $\delta(H) \geq d$.
\end{restatable}

Our proof shows that $\fmain(2,d)$ can be taken to be at most triple exponential in $d$, which is certainly not tight. Observe that, in the case of $G$ being a complete graph, Theorem~\ref{thm:main} implies Ramsey's theorem, hence showing that $\fmain(2,d)$ must be at least as large as the diagonal Ramsey number $\Ram(d+1,d+1)$, which is at least $2^{d/2}$~\cite{erdosBAMS53}. 

\medskip

Combined with a recent result of Kwan {\it et al.}~\cite{KwanComb40} (see Theorem~\ref{thm:KLST} below), Theorem~\ref{thm:main} implies the following for oriented graphs. An {\it oriented graph} is a digraph with at most one arc between any pair of vertices. It is {\it antidirected} if all its vertices are either sources or sinks. The minimum degree of an oriented graph is the minimum degree of its underlying graph. The {\it transitive tournament} of order $r$, denoted by $TT_r$, is the unique acyclic orientation of $K_r$.

\begin{restatable}{theorem}{mainthmoriented}
    \label{thm:mainoriented}
    There exists $\fmainoriented\colon \mathbb{N}^2\to \mathbb{N}$ such that, for all $r,d\in \mathbb{N}$, every oriented graph $D$ with $\delta(D) \geq \fmainoriented(r,d)$ contains $TT_r$ or an induced antidirected subdigraph $H$ with $\delta(H) \geq d$.
\end{restatable}

Using Theorems~\ref{thm:main} and~\ref{thm:mainoriented}, we show that certain extensions of degree-bounded graph classes preserve degree-boundedness. We begin with some definitions.

\medskip

Let $G$ be a graph and $\mathcal{G}$ be a hereditary class of graphs ({\it i.e.} $\mathcal{G}$ is closed under taking induced subgraphs).
We denote by $\tau(G)$ the {\it biclique number} of $G$, that is, the largest integer $s$ such that $G$ contains a copy of $K_{s,s}$ as a subgraph.
The {\it degeneracy} of $G$, denoted by $\delta^\star(G)$, is the least integer $d$ such that every nonempty subgraph $H$ of $G$ contains a vertex of degree at most $d$. 
Observe that $\delta^\star(G)\geq \tau(G)$ trivially holds. A class of graphs is {\it degree-bounded} if, conversely, $\tau$ can be bounded in terms of $\delta^\star$.
That is, $\mathcal{G}$ is degree-bounded if there exists a function $f_{\cal G}\colon \mathbb{N}\to \mathbb{N}$, called a {\it degree-bounding function} of $\mathcal{G}$, such that every graph $H\in \mathcal{G}$ satisfies $\delta^\star(H) \leq f_{\cal G}(\tau(H))$.
Degree-boundedness is a recent topic that attracted a lot of attention in the past few years, see for instance~\cite{mccartySIDMA35,scottJGT102,bonamyJCTB152,bourneufAC24,duJCTB173,giraoIMRN2025}. The interested reader is also referred to~\cite{duEJC104092} for a survey. 

\medskip

Let $k\geq 1$ be an integer. We say that $G$ is {\it $(\mathcal{G},k)$-colourable} if there exists a $k$-colouring $\phi\colon E(G) \to \{1,\dots,k\}$ such that every monochromatic induced subgraph of $G$ belongs to $\mathcal{G}$. 
We define the {\it $k$-extension} $\mathcal{G}^{(k)}$ of $\mathcal{G}$ as the class of all $(\mathcal{G},k)$-colourable graphs. 
The following is a direct consequence of Theorem~\ref{thm:main}.

\begin{corollary}
    \label{cor:main1}
    Let $\mathcal{G}$ be a hereditary class of graphs and $k\geq 1$ be an integer. If $\mathcal{G}$ is degree-bounded, then so is $\mathcal{G}^{(k)}$.
\end{corollary}
\begin{proof}
    If $f_{\mathcal{G}}$ is a degree-bounding function for $\mathcal{G}$, then $x \mapsto \fmain(k,f_{\mathcal{G}}(x))$ is a degree-bounding function for $\mathcal{G}^{(k)}$.
\end{proof}

Similarly, we say that $G$ is {\it $\mathcal{G}$-orientable} if there exists an orientation $D$ of $G$ such that, for every induced antidirected subdigraph $H$ of $D$, the underlying graph of $H$ belongs to $\mathcal{G}$. We define the {\it antidirected extension} $\vec{\mathcal{G}}$ of $\mathcal{G}$ as the class of all $\mathcal{G}$-orientable graphs.
The following is a direct consequence of Theorem~\ref{thm:mainoriented}.

\begin{corollary}
    \label{cor:main2}
    Let $\mathcal{G}$ be a hereditary class of graphs. If $\mathcal{G}$ is degree-bounded, then so is $\vec{\mathcal{G}}$.
\end{corollary}
\begin{proof}
    If $f_{\mathcal{G}}$ is a degree-bounding function for $\mathcal{G}$, then $x \mapsto \fmainoriented(2x,f_{\mathcal{G}}(x))$ is a degree-bounding function for $\vec{\mathcal{G}}$.
\end{proof}

After recalling some notation and preliminary results in Section~\ref{sec:preliminaries}, we discuss applications of Corollaries~\ref{cor:main1} and~\ref{cor:main2} in Section~\ref{sec:applications}, and propose a few further research directions along the way. Proofs of Theorems~\ref{thm:main} and~\ref{thm:mainoriented} are postponed to Section~\ref{sec:proofs}.

\section{Preliminaries}
\label{sec:preliminaries}

\subsection{General notation}

Let $D$ be a digraph. A pair of arcs in opposite directions between the same pair of vertices of $D$ is called a {\it digon}. A {\it simple arc} is an arc that does not belong to any digon.
A digraph without any digon is an {\it oriented graph}. A digraph without any simple arc is {\it symmetric}.
The {\it underlying graph} of $D$ is the undirected graph whose vertex set is $V(D)$ and contains an edge $uv$ whenever there is at least one arc between $u$ and $v$ in $D$. An {\it orientation} of a graph $G$ is any oriented graph whose underlying graph is $G$.

For every integer $\ell \geq 2$, we denote by $\AP_\ell$ the {\it antidirected path} on $\ell$ vertices, that is, the unique antidirected digraph whose underlying graph is $P_\ell$, the path on $\ell$ vertices, and whose initial vertex is a source. We denote by $\vec{P}_\ell$ the {\it directed path} on $\ell$ vertices, that is, the unique orientation of $P_\ell$ with exactly one source and one sink. 
If $\ell\geq 4$ is even, we denote by $\AC_\ell$ the antidirected cycle on $\ell$ vertices, by $\AC_{\geq \ell}$ the set of antidirected cycles on at least $\ell$ vertices, and by $\AC$ the set of all antidirected cycles.

Let $G$ be a graph and $\mathcal{F}$ be any (possibly infinite) set of graphs.
We say that $G$ is {\it $\mathcal{F}$-free} if $G$ does not contain any member of $\mathcal{F}$ as an induced subgraph.
Similarly, for any set $\mathcal{F}$ of oriented graphs, we say that a graph $G$ is $\mathcal{F}$-free if some orientation of $G$ does not contain any member of $\mathcal{F}$ as an induced subdigraph.
With a slight abuse of notation, for an (oriented) graph $F$, we say that $G$ is {\it $F$-free} if it is $\{F\}$-free.
We further say that a graph is {\it triangle-free} if it is $K_3$-free.
A {\it hole} in $G$ is an induced cycle of $G$ of length at least $4$. It is {\it even} if it has even length. We say that a graph is {\it even-hole-free} if it does not contain any even hole, that is, it is $\mathcal{C}_{\rm even}$-free, where $\mathcal{C}_{\rm even}$ is the set of cycles of even length.

\subsection{Basic properties and probabilistic tools}

The {\it average degree} $\Ad(G)$ of a graph $G$ is exactly $2|E(G)|/|V(G)|$.
We omit the proof of the following well-known property, see for instance~\cite[Proposition~1.2.2]{diestel2017}.
\begin{proposition}
    \label{prop:ad_to_mindeg}
    Every graph $G$ with $\Ad(G) \geq d$ contains an induced subgraph $G'$ with $\delta(G') \geq \frac{1}{2}d$. 
\end{proposition}

Our proof of Theorem~\ref{thm:main} contains probabilistic arguments. We make use of the following two well-known bounds.

\begin{proposition}[{\sc Chernoff}]
    \label{prop:chernoff}
    Let $X$ be a random variable following a binomial law with parameters $p \in [0,1]$ and $n \geq 0$, with expectation $\EE(X) = \mu = np$. Then
    \begin{itemize}
        \item $\displaystyle \PP(X \geq (1+\epsilon)\mu) \leq \exp \left(-\epsilon^2 \mu/3\right)$ for any $\epsilon \geq 0$, and
        \item $\displaystyle \PP(X \leq (1-\epsilon)\mu) \leq \exp \left(-\epsilon^2\mu/2\right)$ for any $0 <\epsilon < 1$.
    \end{itemize}
\end{proposition}

\begin{proposition}[{\sc Markov}]
    \label{prop:markov}
    Let $X$ be a nonnegative random variable with $\EE(X)>0$. Then for any $c >0$, $\PP(X\geq c\EE(X))\leq 1/c$.
\end{proposition}

\subsection{Degree- and \texorpdfstring{$\chi$}{chi}-boundedness}

In 2004, Kühn and Osthus~\cite{kuhnComb24} obtained the following seminal result. A {\it subdivision} of a graph $G$ is any graph obtained from $G$ by replacing each edge with a path. It is an {\it even subdivision} if each path has an even number of edges. 

\begin{theorem}[Kühn and Osthus~\cite{kuhnComb24}]
    \label{thm:KO}
    For any fixed graph $H$, the class of graphs that do not contain any induced even subdivision of $H$ is degree-bounded.
\end{theorem}

The fact that, in Theorem~\ref{thm:KO}, the induced subdivision of $K_r$ can be chosen to be even is not explicitly stated in~\cite{kuhnComb24}. To see this, note that it is proved~\cite[Theorem~1.2]{kuhnComb24} that every graph with sufficiently large minimum degree contains $K_{s,s}$ or an induced $1$-subdivision of a graph $H$ with large minimum degree, where a $1$-subdivision is a subdivision in which all edges are subdivided exactly once. Combined with a result of Mader~\cite{maderMA174} (see also~\cite[Theorem~3.5.1]{diestel2000}) stating that every graph with sufficiently large minimum degree contains a subdivision of $K_r$, we indeed obtain Theorem~\ref{thm:KO}.

The degree-bounding function obtained from the original proof of Kühn and Osthus is triple exponential. A natural question is whether this can be strengthened into a polynomial one. Recently, Girão and Hunter~\cite{giraoIMRN2025} settled the problem by showing the following surprising result.

\begin{theorem}[Girão and Hunter~\cite{giraoIMRN2025}]
    \label{thm:girao_hunter}
    Any hereditary degree-bounded class $\mathcal{G}$ has a degree-bounding function that is a polynomial.
\end{theorem}

Even though degree-boundedness might present some similarities with the so-called notion of {\it $\chi$-boundedness}, Theorem~\ref{thm:girao_hunter} makes an important distinction between them. We say that $\mathcal{G}$ is {\it $\chi$-bounded} if there exists a function $f_{\cal G}\colon \mathbb{N}\to \mathbb{N}$, called a {\it $\chi$-bounding function} of $\mathcal{G}$, such that any graph $G\in \mathcal{G}$ satisfies $\chi(G) \leq f_{\cal G}(\omega(G))$.
The analogue of Theorem~\ref{thm:girao_hunter} for $\chi$-boundedness was indeed a famous conjecture posed by Esperet~\cite{esperetHdr} and disproved recently by Bria{\'n}ski, Davies, and Walczak~\cite{brianskiComb44}. We say that $\mathcal{G}$ is {\it polynomially} (respectively {\it linearly}) $\chi$-bounded if $\mathcal{G}$ admits a polynomial (respectively linear) $\chi$-bounding function.

\medskip

Observe that Theorem~\ref{thm:KO}, applied with $H=K_3$, directly implies that the class of even-hole-free graphs is degree-bounded. Noticing that any induced copy of $K_{2,2}$ is an even hole, it even shows that even-hole-free graphs $G$ satisfy $\delta(G) \leq f(\omega(G))$ for some function $f$. In particular, since the class of even-hole-free graphs is hereditary, it follows that it is $\chi$-bounded. Furthermore, using Theorem~\ref{thm:girao_hunter}, it even shows that it is polynomially $\chi$-bounded.
Actually, Chudnovsky and Seymour~\cite{addarioJCTB98,chudnovskyJCTB161} proved the following more precise result.

\begin{theorem}[\cite{addarioJCTB98,chudnovskyJCTB161}]
    \label{thm:EH_free}
    Every even-hole-free graph $G$ satisfies $\delta(G)\leq 2\omega(G)-2$.
\end{theorem}

\subsection{Extremal combinatorics}

For integers $a,b\in \mathbb{N}$, we denote by $\Ram(a,b)$ the least integer $r$ such that every graph of order $r$ contains a clique of size $a$ or an independent set of size $b$.

\begin{theorem}[{\sc Ramsey}~\cite{ramsey1930}]
    \label{thm:ramsey}
    For all integers $a,b$, $\Ram(a,b)$ exists and $\Ram(a,b)\leq \binom{a+b-2}{a-1}$.
\end{theorem}

The upper bound on $\Ram(a,b)$ above follows from the celebrated proof of Erd\H{o}s and Szekeres, see for instance~\cite[Theorem~12.5]{bondy2008}. 

The following is a directed analogue of Ramsey's Theorem. A {\it tournament} is any orientation of a complete graph. Recall that the {\it transitive tournament} of order $r$, denoted by $TT_r$, is the unique acyclic tournament of order $r$.

\begin{theorem}[Erd\H{o}s and Moser~\cite{erdos1964}]
    \label{thm:erdos_moser}
    Every tournament of order $2^{r-1}$ contains $TT_r$.
\end{theorem}

We further need the following celebrated result due to K\H{o}v\'ari, S\'os, and Tur\'an.

\begin{theorem}[{K\H{o}v\'ari, S\'os, and Tur\'an}~\cite{kovariCM3}]
    \label{thm:kst}
    For every $s\in \mathbb{N}$, there exists $c_s$ such that every graph of order $n$ without any copy of $K_{s,s}$ has at most $c_s\cdot n^{2-\frac{1}{s}}$ edges.
\end{theorem}

The following is a direct consequence of Theorem~\ref{thm:kst}. 

\begin{corollary}
    \label{cor:kst}
    For all integers $s,k\geq 1$, there exists an integer $t$ such that every $k$-edge-colouring of $K_{t,t}$ yields a monochromatic copy of 
    $K_{s,s}$.
\end{corollary}
\begin{proof}
    Consider any $k$-edge-colouring of $K_{t,t}$ with colour classes $E_1,\dots,E_k$. 
    Then, for some $i\in [k]$, $E_i$ has size at least $\frac{1}{k}{t}^{2}$ edges. 
    Remove all the other edges. What remains is a graph $G$ of order $2t$ with at least $\frac{1}{k}{t}^{2}$ edges. If $t$ is large enough compared to $k$ and $s$, Theorem~\ref{thm:kst} guarantees that the underlying graph of $G$ contains a copy of $K_{s,s}$. 
\end{proof}

We finally need the following result due to Kwan, Letzter, Sudakov, and Tran~\cite{KwanComb40}, which states that, in any graph with large minimum degree, one can always extract a large clique or an induced bipartite graph with large minimum degree, thereby addressing a conjecture of Esperet, Kang, and Thomassé~\cite{esperetCPC28}.

\begin{theorem}[Kwan {\it et al.}~\cite{KwanComb40}]
    \label{thm:KLST}
    There exists a function $\fKLST : \mathbb{N}^2\to \mathbb{N}$ such that, for all $r,d\in \mathbb{N}$, every graph $G$ with $\delta(G)\geq \fKLST(r,d)$ contains $K_r$ or an induced bipartite graph $H$ with $\delta(H) \geq d$.
\end{theorem}

\section{Applications to degree- and \texorpdfstring{$\chi$}{chi}-boundedness}
\label{sec:applications}

In this section, we illustrate how Corollaries~\ref{cor:main1} and~\ref{cor:main2} can be used to show that some well-known graph classes are degree- or (polynomially) $\chi$-bounded. 

\subsection{A directed version of Kühn and Osthus' theorem}

We first notice that Corollary~\ref{cor:main2} can be combined with Theorem~\ref{thm:KO} to obtain the following directed analogue of Kühn and Osthus' theorem.
Let $H$ be an undirected graph.
An {\it antidirected subdivision} of $H$ is any antidirected digraph $D$ whose underlying graph is an even subdivision of $H$.
Observe that, by definition of an even subdivision, the branch vertices of $D$ are either all sources or all sinks.
Given an integer $s$, we denote by $\directedKss$ the unique antidirected digraph whose underlying graph is~$K_{s,s}$. 

\begin{corollary}
    \label{cor:substructure:oriented}
    There exists $\fcorsubstructureoriented \colon \mathbb{N}^2\to \mathbb{N}$ such that the following holds.
    Every oriented graph $D$ with $\delta(D) \geq \fcorsubstructureoriented(r,s)$ contains a copy of $\directedKss$ or an induced antidirected subdivision of $K_r$.
\end{corollary}
\begin{proof}
    Let $\mathcal{S}_r$ denote the class of graphs without any induced even subdivision of $K_r$. By Theorem~\ref{thm:KO}, $\mathcal{S}_r$ is degree-bounded, and so is $\vec{\mathcal{S}_r}$ by Corollary~\ref{cor:main2}.

    Now, let $D$ be an oriented graph with minimum degree sufficiently large compared to $r$ and $s$, and assume that $D$ does not contain an induced antidirected subdivision of $K_r$. Then, by definition, $D$ belongs to $\vec{\mathcal{S}_r}$, implying that the underlying graph of $D$ contains a copy $K$ of $K_{t,t}$ for some $t$ large enough with respect to $s$.
    Let $(A,B)$ denote the bipartition of $K$. 
    It follows from Corollary~\ref{cor:kst} that $D$ contains $\directedKss$ by colouring the edges of $K$ with $1$ if they go from $A$ to $B$ and with $2$ otherwise.
\end{proof}

Observe that the analogue of Corollary~\ref{cor:substructure:oriented} does not hold when we look for a {\it directed subdivision} of $K_r$ -- that is, a digraph obtained from $K_r$ by replacing each edge with a directed path -- instead of an antidirected one. To see this, it is sufficient to take any bipartite graph $(A,B)$ with arbitrarily large girth and minimum degree, and to orient all the edges from $A$ to $B$. In the obtained oriented graph, every subdigraph is antidirected, while every directed subdivision of $K_3$ contains a vertex which is neither a source nor a sink.

We actually believe that a more precise version of Corollary~\ref{cor:substructure:oriented} holds. A {\it restricted antidirected subdivision} of $H$ is an antidirected subdivision of $H$ in which all branch vertices are sources. We pose the following conjecture.

\begin{conjecture}
    Every oriented graph $D$ with $\delta(D) \geq \fmainoriented(2s,\fKO (r,s))$ contains a copy of $\directedKss$ or an induced restricted antidirected subdivision of $K_r$.
\end{conjecture}

\subsection{Odd-signable graphs}

A class of graphs containing even-hole-free graphs is  the {\it odd-signable} graphs. A graph $G$ is odd-signable if its edges can be labelled with $\{0,1\}$ in such a way that, for every induced cycle, the sum of its labels is odd.
This class was mainly introduced for its structural similarities with even-hole-free graphs, and was used to show that even-hole-free graphs can be recognised in polynomial time, see for instance~\cite{confortiJGT30,confortiJGT39,vuvskovicAADM4,abrishamiJCTB157}.

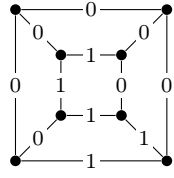
\begin{figure}
\centering
\begin{tikzpicture}[scale=1]
    \node[vertex] (000) at (0,0) {};
    \node[vertex] (001) at (2,0) {};
    \node[vertex] (010) at (0,2) {};
    \node[vertex] (011) at (2,2) {};
    \node[vertex] (100) at (0.6,0.6) {};
    \node[vertex] (101) at (1.4,0.6) {};
    \node[vertex] (110) at (0.6,1.4) {};
    \node[vertex] (111) at (1.4,1.4) {};

    \draw (000) -- node[midway, fill=white, inner sep=1pt]  {\scriptsize 1} (001);
    \draw (000) -- node[midway, fill=white, inner sep=2pt] {\scriptsize 0} (010);
    \draw (000) -- node[midway, fill=white, inner sep=1pt] {\scriptsize 0} (100);
    \draw (001) -- node[midway, fill=white, inner sep=1pt] {\scriptsize 1} (101);
    \draw (001) -- node[midway, fill=white, inner sep=2pt] {\scriptsize 0} (011);
    \draw (010) -- node[midway, fill=white, inner sep=1pt] {\scriptsize 0} (110);
    \draw (010) -- node[midway, fill=white, inner sep=1pt] {\scriptsize 0} (011);
    \draw (011) -- node[midway, fill=white, inner sep=1pt] {\scriptsize 0} (111);
    \draw (100) -- node[midway, fill=white, inner sep=1pt] {\scriptsize 1} (101);
    \draw (100) -- node[midway, fill=white, inner sep=2pt] {\scriptsize 1} (110);
    \draw (101) -- node[midway, fill=white, inner sep=2pt] {\scriptsize 0} (111);
    \draw (110) -- node[midway, fill=white, inner sep=1pt] {\scriptsize 1} (111);
\end{tikzpicture}
\caption{An edge-labelling of the cube witnessing that it is odd-signable.}
\label{fig:odd_signable:cube}
\end{figure}

Note that the result of Chudnovsky and Seymour (Theorem~\ref{thm:EH_free}) does not extend to odd-signable graphs, as there exist odd-signable graphs without any bisimplicial vertex (for instance, the cube, see Figure~\ref{fig:odd_signable:cube}).
Still, by definition, odd-signable graphs are in $\mathcal{EH}^{(2)}$, where $\mathcal{EH}$ denotes the class of even-hole-free graphs.
Therefore, combining Corollary~\ref{cor:main1} with Theorem~\ref{thm:EH_free}, we obtain that odd-signable graphs are degree-bounded. It turns out that they also are polynomially $\chi$-bounded.

\begin{corollary}
    \label{cor:chiboundedness_general_evenhole}
    For every $k\in \mathbb{N}$, $\Geh{k}$ is degree-bounded and polynomially $\chi$-bounded.
\end{corollary}
\begin{proof}
    Let us fix $k\in \mathbb{N}$. Clearly, $\Geh{k}$ is degree-bounded by Corollary~\ref{cor:main1} and Theorem~\ref{thm:EH_free}, so we only have to show polynomial $\chi$-boundedness. 
    
    Let $g_k\colon \mathbb{N}\to \mathbb{N}$ be a polynomial degree-bounding function of $\Geh{k}$, whose existence is guaranteed by Theorem~\ref{thm:girao_hunter}. Let $s_k$ be the least integer such that any $k$-edge-colouring of $K_{s_k,s_k}$ contains a monochromatic copy of $K_{2,2}$, whose existence is guaranteed by Corollary~\ref{cor:kst}.
    
    Let $G$ be any graph in $\Geh{k}$. We claim that it satisfies $\delta(G) < g_k(\Ram(\omega(G)+1,s_k))$. Indeed, if $\delta(G) \geq g_k(\Ram(\omega(G)+1,s_k))$, then, by definition of $g_k$, $\tau(G)\geq \Ram(\omega(G)+1,s_k)$. It follows that $G$ contains an induced copy of $K_{s_k,s_k}$. By definition of $s_k$, $G$ contains a monochromatic induced copy of $K_{2,2}$, a contradiction.

    We just proved that every graph $G\in \Geh{k}$ satisfies $\delta(G) \leq  g_k(\Ram(\omega(G)+1,s_k))-1$. It follows by induction ($\Geh{k}$ being hereditary) that every such graph $G$ satisfies $\chi(G) \leq g_k(\Ram(\omega(G)+1,s_k))$. Note that, by Theorem~\ref{thm:ramsey}, when $k$ is fixed, $\Ram(\omega+1,s_k)$ is polynomial in $\omega$, and so is $g_k(\Ram(\omega+1,s_k))$. The result follows.
\end{proof}

Even though it is polynomial, the $\chi$-bounding function we obtain is certainly not optimal. We ask whether a linear one exists.

\begin{problem}
    Is $\Geh{k}$ linearly $\chi$-bounded for every $k$?
\end{problem}

\subsection{Analogues of Gy\'arf\'as-Sumner's Conjecture}

Let $F$ be any graph such that the class of $F$-free graphs is $\chi$-bounded. Since there exist graphs of arbitrarily large girth and chromatic number~\cite{erdosCJM11}, $F$ must be a forest. The celebrated Gy\'arf\'as-Sumner's Conjecture~\cite{gyarfasIFS2,sumner1981} states that the converse is also true. 

\begin{conjecture}[{\sc Gy\'arf\'as-Sumner}~\cite{gyarfasIFS2,sumner1981}]
    For every graph $F$, the class of $F$-free graphs is $\chi$-bounded if and only if $F$ is a forest.
\end{conjecture}

For degree-boundedness, the analogous problem has been settled by Kierstead and Penrice~\cite{kiersteadJGT18} (even though it is not explicitly stated in the paper, see~\cite[Theorem~3.15]{duEJC104092}).

\begin{theorem}[Kierstead and Penrice~\cite{kiersteadJGT18}]
    \label{thm:kierstead}
    Let $F$ be a graph.
    The class of $F$-free graphs is degree-bounded if and only if $F$ is a forest.
\end{theorem}

We derive the analogous result for oriented graphs. By an {\it antidirected forest}, we mean any antidirected digraph whose underlying graph is a forest.

\begin{corollary}
    \label{cor:oriented_degree_gyarfassumner}
    Let $F$ be an oriented graph. The class of $F$-free graphs is degree-bounded if and only if $F$ is an antidirected forest.
\end{corollary}
\begin{proof}
    If, for some oriented graph $F$, the class of $F$-free graphs is degree-bounded, then $F$ must be an antidirected forest, as there exist antidirected digraphs with arbitrarily large girth and minimum degree. To see this, take any bipartite graph with arbitrarily large girth and minimum degree, and orient it from one part of the bipartition to the other. 

    Let us now fix an antidirected forest $F$ with underlying graph $F'$. We may assume that $F'$ is connected, for otherwise we choose one arbitrary sink in every connected component of $F$, and add a new vertex dominating all these sinks.
    Let $v$ be an arbitrary vertex of $F$, and let $F^\star$ be the forest obtained from two disjoint copies of $F'$ by adding an edge between the two copies of $v$. 
    Observe that every antidirected digraph whose underlying graph is $F^\star$ contains an induced copy of $F$.
    Therefore, the class of $F$-free graphs is a subclass of the antidirected extension of $F^\star$-free graphs, and the result follows from the combination of Corollary~\ref{cor:main2} and Theorem~\ref{thm:kierstead} (applied to~$F^\star$).
\end{proof}

We can further derive the analogous result for digraphs. A {\it symmetric forest} is any symmetric digraph whose underlying graph is a forest. 
Let $\mathcal{F}$ be any set of digraphs. We say that a graph $G$ is {\it $\mathcal{F}$-difree} if there exists a digraph $D$ whose underlying graph is $G$ and $D$ does not contain any induced copy of any member of $\mathcal{F}$.
Observe that the following is a direct consequence of Theorem~\ref{thm:main}.

\begin{corollary}
    \label{cor:homogeneous_digraph}
    Every digraph $D$ with $\delta(D) \geq \fmain(2,d)$ contains an induced subdigraph $H$ with $\delta(H) \geq d$ that is either a symmetric digraph or an oriented graph.
\end{corollary}
\begin{proof}
    Let $G$ be the underlying graph of $D$ and $\phi$ be the $2$-edge-colouring of $G$ where $\phi(\{u,v\}) = 1$ if $\{uv,vu\}$ is a digon of $D$, and $\phi(\{u,v\}) = 2$ otherwise. The result now follows from Theorem~\ref{thm:main} applied to $G$.
\end{proof}

Therefore, combining Corollaries~\ref{cor:homogeneous_digraph} and~\ref{cor:oriented_degree_gyarfassumner} with Theorem~\ref{thm:kierstead}, we obtain the following. 

\begin{corollary}
    For every symmetric forest $F_1$ and antidirected forest $F_2$, the class of $\{F_1,F_2\}$-difree graphs is degree-bounded.
\end{corollary}

Again, note that this is optimal in the following sense: if $\mathcal{F}$ is a finite set of digraphs such that the class of $\mathcal{F}$-difree graphs is degree-bounded, then $\mathcal{F}$ contains at least one symmetric forest and at least one antidirected forest. 

Regarding $\chi$-boundedness, it is not clear what are the oriented graphs $F$ such that $F$-free graphs are $\chi$-bounded. If $F$ is such an oriented graph, then its underlying graph is a forest, but the converse is not true. As noticed by Gy\'arf\'as~\cite{gyarfasDM79}, the class of $\AP_4$-free graphs contains the shift graphs, which are known to be triangle-free graphs with unbounded chromatic number~\cite{erdosTG68}.
Similarly, Kierstead and Trotter~\cite{kiersteadDM101} gave a construction, similar to the one of Zykov~\cite{zykov1952}, of $\vec{P}_4$-free triangle-free graphs with unbounded chromatic number.

On the positive side, Chudnovsky and Seymour~\cite{chudnovskyEJC76} proved that excluding
\tikz{
\node[vertex,white] (u) at (0,-0.05) {};
\node[smallvertex] (a) at (0,0) {};
\node[smallvertex] (b) at (0.5,0){};
\node[smallvertex] (c) at (1,0){};
\node[smallvertex] (d) at (1.5,0){};
\draw[arc](a) to (b);
\draw[arc](c) to (b);
\draw[arc](d) to (c);
}, its directional dual 
\tikz{
\node[vertex,white] (u) at (0,-0.05) {};
\node[smallvertex] (a) at (0,0) {};
\node[smallvertex] (b) at (0.5,0){};
\node[smallvertex] (c) at (1,0){};
\node[smallvertex] (d) at (1.5,0){};
\draw[arc](b) to (a);
\draw[arc](b) to (c);
\draw[arc](c) to (d);
}, or any orientation of a star yields a $\chi$-bounded class of graphs, thereby confirming a conjecture of Aboulker, Bang-Jensen, Bousquet, Charbit, Havet, Maffray, and Zamora~\cite{aboulkerJGT89}. The following problem is open, and the interested reader is referred to the last section of~\cite{aboulkerJGT89} for further details.

\begin{problem}
    Characterise the oriented forests $F$ such that $F$-free graphs are $\chi$-bounded.
\end{problem}

\subsection{Antidirected cycles and Burling graphs}

Observe that $\AC$-free graphs form another generalisation of even-hole-free graphs. As a consequence of Corollary~\ref{cor:substructure:oriented}, we obtain the following.

\begin{corollary}
    \label{cor:degree_boundedness_antidirected_cycles}
    For every integer $\ell$, the class of $\AC_{\geq \ell}$-free graphs is degree-bounded.
\end{corollary}
\begin{proof}
    Let $G$ be any $\AC_{\geq \ell}$-free graph, and let $\tau(G) = \tau$.
    We claim that 
    \[ \delta(G) < \fmainoriented\Big(2\tau+2,\fKO(\ell,\tau+1)\Big),\] hence implying the result.
    Assume that this is not the case, and let $D$ be any orientation of $G$ without any induced antidirected cycle of length at least $\ell$. By Corollary~\ref{cor:substructure:oriented}, $D$ contains an induced antidirected subdivision of $K_\ell$, which we denote $K$.

    Let us label $u_1,\dots,u_\ell$ the vertices of $K_\ell$, and for any $i<j$, let us denote by $P_{i,j}$ the antidirected path connecting $u_i$ and $u_j$ in $K$. The concatenation of $P_{1,2},P_{2,3},\dots,P_{\ell-1,\ell},P_{\ell,1}$ is an induced antidirected cycle of length at least $2\ell$ in $D$, a contradiction.
\end{proof}

By noticing that $\directedK_{2,2}$ is exactly $\AC_4$, and that all orientations of $K_{s,s}$ for some $s\in \mathbb{N}$ contain a copy of $\directedK_{2,2}$, we derive the following.

\begin{corollary}
    \label{cor:AC}
    The class of $\AC$-free graphs is polynomially $\chi$-bounded.
\end{corollary}
\begin{proof}
    By Corollary~\ref{cor:degree_boundedness_antidirected_cycles}, $\AC$-free graphs are degree-bounded. By Theorem~\ref{thm:girao_hunter}, let $g$ be any polynomial degree-bounding function of $\AC$-free graphs.
    
    Let $s$ be the least integer such that any $2$-edge-colouring of $K_{s,s}$ yields a monochromatic copy of $K_{2,2}$, the existence of which is guaranteed by Corollary~\ref{cor:kst}. Equivalently, any orientation of $K_{s,s}$ contains a copy of $\AC_4$. To see this, denote by $(A,B)$ the bipartition of $K_{s,s}$. To any orientation $D$ of $K_{s,s}$ we can associate the $2$-edge-colouring $\phi_D$ in which the edges going from $A$ to $B$ are coloured~$1$, and the ones from $B$ to $A$ are coloured $2$. Now $\phi_D$ yields a monochromatic copy of $K_{2,2}$ if and only if $D$ contains a copy of $\directedK_{2,2}$. 

    Let us fix any $\AC$-free graph $G$ with clique number $\omega$. We claim that $\delta(G) <  g(\Ram(\omega+1,s))$, and it follows that $\chi(G) \leq g(\Ram(\omega+1,s))$ (the class of $\AC$-free graphs being hereditary). 
    Note that, $s$ being a constant, $\Ram(\omega+1,s)$ is polynomial in $\omega$, and so is $g(\Ram(\omega+1,s))$.
    
    Suppose for a contradiction that $\delta(G) \geq  g(\Ram(\omega+1,s))$. It follows that $\tau(G) \geq \Ram(\omega+1,s)$. Since $G$ has clique number $\omega$, $G$ contains an induced copy of $K_{s,s}$. By choice of $s$, any orientation of $G$ contains an induced copy of $\AC_4$, a contradiction. The result follows.
\end{proof}

Finally, we notice that Corollary~\ref{cor:degree_boundedness_antidirected_cycles} can be combined with a recent result of Pournajafi and Trotignon~\cite{pournajafiEJC116} on {\it Burling graphs} (as defined in~\cite{pournajafiEJC116}). Burling graphs form a hereditary class of triangle-free graphs with unbounded chromatic number, due to Burling~\cite{burlingthesis}.
These graphs are particular intersection graphs of various geometrical objects, such as axis-aligned boxes in $\mathbb{R}^3$~\cite{burlingthesis}, line segments in the plane~\cite{pawlikJCTB105}, or frames~\cite{pawlikDCG50,KrawczykDCG53,chalopinEJC23}. Burling graphs attracted a lot of attention since Pawlik, Kozik, Krawczyk, Lason, Micek, Trotter, and Walczak~\cite{pawlikJCTB105} proved that they form a counterexample to a celebrated conjecture of Scott~\cite{scottJGT24}.

Pournajafi and Trotignon~\cite{pournajafiEJC116} recently gave five new characterizations of Burling graphs, three of them being geometrical. With a fourth one, called {\it derived graphs}, they obtained the following (see also~\cite[Corollary~7.31]{pournajafithesis}). 

\begin{theorem}[Pournajafi and Trotignon~\cite{pournajafiEJC116}]
    Every Burling graph admits an orientation in which each hole has exactly two sources.
\end{theorem}

In particular, the result above implies that every Burling graph admits an orientation in which all induced antidirected cycles have length exactly $4$. 
In other words, the class of $\AC_{\geq 6}$-free graphs contains Burling graphs. It follows that  $\AC_{\geq 6}$-free graphs are not $\chi$-bounded and, together with Corollary~\ref{cor:degree_boundedness_antidirected_cycles}, that Burling graphs are degree-bounded.

\begin{corollary}
    \label{cor:burling}
    Burling graphs are degree-bounded.
\end{corollary}

Corollary~\ref{cor:burling} was already proved (see~\cite[Theorem~3.8]{duEJC104092}) with a geometrical approach. Indeed, as mentioned above, Pawlik, Kozik, Krawczyk, Laso\'n, Micek, Trotter, and Walczak~\cite{pawlikJCTB105} proved that Burling graphs can be seen as intersection graphs of line segments in the plane. Fox and Pach~\cite{foxCPC19} proved that the more general class of {\it String graphs} (that is, intersection graphs of continuous arcs in the plane) is degree-bounded. These two results together show that Burling graphs are degree-bounded.
Observe that Corollary~\ref{cor:burling} is in contrast with a recent result of Davies~\cite{daviesAC2021}, who proved the existence of intersection graphs of axis-aligned boxes in $\mathbb{R}^3$ with arbitrarily large girth and chromatic number, whereas Burling graphs of girth $5$ have bounded degeneracy.

\section{Dense monochromatic induced subgraphs}
\label{sec:proofs}

This section is dedicated to the proofs of Theorems~\ref{thm:main} and~\ref{thm:mainoriented}.
We need the following intermediate result due to Kühn and Osthus, which is shown in the proof of~\cite[Lemma~10]{kuhnComb24}, even though it is not explicitly stated this way. We give the proof for completeness.

\begin{lemma}[Kühn and Osthus~\cite{kuhnComb24}]
    \label{lemma:KO1}
    Let $\Gamma,d$ be real numbers with $\Gamma> 16d\geq 32$. Every bipartite graph $G$ with $\Ad(G) =\Gamma$ contains an induced subgraph $G^\star$ with bipartition $(A^\star, B^\star)$ such that $|A^\star| \geq \frac{\Gamma}{128d}|B^\star|$ and, for every $a\in A^\star$, $4d \leq d_{G^\star}(a) \leq 64d$.
\end{lemma}
\begin{proof}
    We assume that $\Ad(G) = \Gamma$ and that $\Ad(H)\leq \Gamma$ for all subgraphs $H$ of $G$, for otherwise we may prove the result on some induced subgraph of $G$ with larger average degree.

    Let $G'$ be an induced subgraph of $G$ with $\delta(G') \geq \Gamma/2$ and bipartition $(A,B)$, the existence of which is guaranteed by Proposition~\ref{prop:ad_to_mindeg}. 
    We assume without loss of generality that $|A| \geq |B|$. In particular, this implies
    \begin{equation}
        \label{eq:KO1}
       \sum_{a\in A}d_{G'}(a) = |E(G')| \leq \frac{\Gamma}{2}(|A|+|B|) \leq \Gamma \cdot |A|.
    \end{equation}
    Let $W$ be the set of vertices $a\in A$ whose degree is larger than $2\Gamma$ in $G'$. Note that $|E(G')| > |W| \cdot 2\Gamma$, which together with~\eqref{eq:KO1} implies $|W| < \frac{1}{2}|A|$.

    Let $B^\star$ be a random subset of $B$ obtained by including each vertex $b\in B$ independently with probability $16d/\Gamma$. Let $A^\star$ be the set of vertices $a\in A \setminus W$ such that $4d \leq |N(a) \cap B^\star| \leq 64d$.
    
    Let us fix any $a\in A\setminus W$ and estimate the probability that $a\in A^\star$.
    For this, let $X_a = |N(a)\cap B^\star|$. Note that $X_a$ follows a binomial law with parameters $16d/\Gamma$ and $d(a) \in [\Gamma/2, 2\Gamma]$. In particular, $8d\leq \EE(X_a)\leq 32d$. Therefore, by Chernoff's Inequality (Proposition~\ref{prop:chernoff}), we have
    \[
        \PP(X_a \geq 64d) \leq \PP(X_a \geq 2\EE(X_a)) \leq e^{-8d/3}
    \]
    and
    \[
    \PP(X_a \leq 4d) \leq \PP(X_a \leq \EE(X_a)/2) \leq e^{-d}.
    \]
    Let $\Tilde{A}$ be the vertices in $A\setminus (A^\star\cup W)$.
    By the union bound, it follows that 
    \begin{equation*}
        \label{eq:KO1:prob_bad}
        \PP(a\in \Tilde{A}) = \PP\Big((X_a<4d) \cup (X_a> 64d)\Big) \leq \frac{1}{4}.
    \end{equation*}
    By linearity of expectation, this implies
    \[
    \EE(|\Tilde{A}|) \leq \frac{1}{4}|A\setminus W|.
    \]
    Hence, by Markov's Inequality (Proposition~\ref{prop:markov}), we obtain
    \begin{equation}
        \label{eq:KO1:sizeAstar}
        \PP\Big(|A^\star| \leq \frac{1}{4}|A|\Big)
        \leq \PP\Big(|A^\star| \leq \frac{1}{2}|A\setminus W| \Big) =
        \PP\Big(|\Tilde{A}| \geq \frac{1}{2}|A\setminus W| \Big)
        \leq \PP\Big(|\Tilde{A}| \geq 2\EE(|\Tilde{A}|) \Big) \leq \frac{1}{2}.
    \end{equation}
    Let us finally estimate the size of $B^\star$. 
    Again, $|B^\star|$ follows a binomial law of parameters $16d/\Gamma$ and $|B|$. We thus have $\EE(|B^\star|) = \frac{16d}{\Gamma}|B| \leq \frac{16d}{\Gamma}|A|$ and $\EE(|B^\star|) \geq 8d$ since $|B|\geq \Gamma/2$ (which holds because every vertex in $G'$ has degree at least $\Gamma/2$). By Chernoff's Inequality (Proposition~\ref{prop:chernoff}), we have
    \begin{equation}
        \label{eq:KO1:sizeBstar}
        \PP\Big(|B^\star| \geq \frac{32d}{\Gamma} |A|\Big) \leq \PP\Big(|B^\star| \geq 2\EE(|B^\star|)\Big) \leq e^{-8d/3} \leq \frac{1}{4}.
    \end{equation}
    By the union bound, it follows from~\eqref{eq:KO1:sizeAstar} and~\eqref{eq:KO1:sizeBstar} that, with positive probability, $|A^\star|\geq \frac{1}{4}|A| \geq \frac{\Gamma}{128d}|B^\star|$. The result follows.
\end{proof}

We are now ready to prove Theorem~\ref{thm:main}, that we first recall here for convenience.

\mainthm*
\begin{proof}
    Note that the existence of $\fmain(k,d)$ is clear when $k\leq 1$. Furthermore, it is sufficient to show the existence of $\fmain(2,d)$ for every $d\in \mathbb{N}$. The existence of $\fmain(k,d)$ for $k\geq 3$ and $d\in \mathbb{N}$ then easily follows by induction on $k$, by taking $\fmain(k,d) = \fmain({k-1},f(2,d))$, and identifying two colour classes in the given $k$-edge-coloured graph $G$.

    Let us fix $d\in \mathbb{N}$. We may assume $d\geq 2$, the existence of $\fmain(2,d)$ being trivial otherwise. We define $g(d) = d\cdot 2^{64d+8}$. We will prove the result for 
    \[
        \fmain(2,d) = \fKLST \Big( \Ram(d+1,d+1), g(d) \Big),
    \]
    Let us fix a graph $G$ with minimum degree at least $\fmain(2,d)$. 
    By Theorem~\ref{thm:KLST}, either $G$ contains a clique of size $\Ram(d+1,d+1)$ or an induced bipartite graph $G'$ with minimum degree at least $g(d)$. In the former case, $G$ contains a monochromatic clique $C$ of size $d+1$, and $G[C]$ is a  monochromatic induced subgraph with minimum degree $d$. We may thus assume the latter.

    Observe that $\Ad(G') \geq \delta(G') \geq g(d) > 16d \geq 32$. Therefore, by Lemma~\ref{lemma:KO1}, $G'$ contains an induced bipartite graph $H$ with bipartition $(A, B)$ such that $|A| \geq 2^{64d+1}|B|$ and every vertex $a\in A$ satisfies 
    \[
    4d \leq d_{H}(a) \leq 64d.
    \]
    For $i\in \{1,2\}$ and every vertex $v$, let $d^i(v)$ denote the number of edges  coloured $i$ incident to $v$. Let $A_1$ be the set of vertices $a\in A$ such that $d_H^1(a) \geq d_H^2(a)$, and similarly let $A_2$ be the set of vertices $a\in A$ such that $d_H^2(a) \geq d_H^1(a)$.
    We assume without loss of generality that $|A_1| \geq |A_2|$, and in particular $|A_1| \geq \frac{1}{2}|A| \geq 2^{64d}|B|$.
    For every $a\in A_1$, we fix a set $B_a$ of exactly $2d$ vertices in $B$, such that every edge $ab$ for $b\in B_a$ is coloured $1$. This is possible by the definition of $A_1$.

    Let $B^\star$ be a subset of $B$ chosen uniformly at random. In other words, each vertex $b\in B$ belongs to $B^\star$ with probability $1/2$. Let $A^\star$ be the set of vertices $a\in A_1$ for which $N_H(a) \cap B^\star = B_a$. Note that, by definition, $G^\star = G[A^\star \cup B^\star]$ is monochromatic, and every vertex in $A^\star$ has degree exactly $2d$ in $G^\star$. 
    Let us now estimate the expected size of $A^\star$. For every $a\in A_1$, we have
    \[
    \PP(a\in A^\star) = \left(\frac{1}{2}\right)^{d_H(a)} \geq \left(\frac{1}{2}\right)^{64d}.
    \]
    Therefore, by linearity of expectation, 
    \[
    \EE(|A^\star|) \geq 2^{-64d}\cdot |A_1| \geq |B| \geq |B^\star|.
    \]
    It follows that, with positive probability, we indeed have $|A^\star|\geq |B^\star|$. For this particular choice, $G^\star$ is a monochromatic induced subgraph of $G$ with average degree
    \[
        \Ad(G^\star) = \frac{2d\cdot |A^\star|}{|A^\star|+|B^\star|} \geq 2d.
    \]
    The result then follows from Proposition~\ref{prop:ad_to_mindeg}.
\end{proof}

We now prove Theorem~\ref{thm:mainoriented}, that we first recall here for convenience.

\mainthmoriented*
\begin{proof}
    Let us fix $r$ and $d$. We prove the result for $\fmainoriented(r,d) = \fKLST(2^{r-1},\fmain(2,d))$. Let us thus fix an orientation $D$ of a graph $G$ with $\delta(G) \geq \fmainoriented(r,d)$.
    
    By Theorem~\ref{thm:KLST}, either $G$ contains a clique of size $2^{r-1}$ or an induced bipartite graph $G'$ with minimum degree at least $\fmain(2,d)$. In the former case, $D$ contains a tournament of order $2^{r-1}$, which by Theorem~\ref{thm:erdos_moser} implies that $D$ contains $TT_r$. We may thus assume the latter.

    Let $D'=D[V(G')]$,  and let $(A,B)$ be the bipartition of $G'$. Let $\phi_D$ be the $2$-edge-colouring of $G'$ where edges oriented from $A$ to $B$ in $D'$ are coloured $1$, and the others are coloured $2$. By definition, every monochromatic subgraph of $G'$ induces an antidirected digraph on $D'$. The result thus follows from~Theorem~\ref{thm:main}.
\end{proof}

\section*{Acknowledgments}

We are thankful to Marthe Bonamy, Clément Rambaud, and Nicolas Trotignon for stimulating discussions and helpful suggestions.

\end{document}